\documentclass[10pt, leqno]{amsart}

\usepackage{amssymb,amsfonts,amscd,amsbsy}
\usepackage[mathscr]{eucal}

\theoremstyle{plain}
\newtheorem{thm}{Theorem}[section]

\newtheorem{pro}[thm]{Proposition}

\theoremstyle{definition}

\newtheorem{rem}[thm]{Remark}

\newtheorem{conj}[thm]{Conjecture}

\begin{document}

\title[Binary Quadratic Forms]{A remark on a conjecture of Borwein and Choi}
\author[Robert Osburn]{Robert Osburn}

\address{Department of Mathematics $\&$ Statistics, Queen's University, Kingston, Ontario, Canada K7L 3N6}

\email{osburnr@mast.queensu.ca}


\begin{abstract}
We prove the remaining case of a conjecture of Borwein and Choi
concerning an estimate on the square of the number of solutions to
$n=x^2+Ny^2$ for a squarefree integer $N$.
\end{abstract}

\maketitle
\section{Introduction}

We consider the positive definite quadratic form $Q(x,y)=x^2+Ny^2$
for a squarefree integer $N$. Let $r_{2,N}(n)$ denote the number
of solutions to $n=Q(x,y)$ (counting signs and order). In this
note, we estimate

\begin{center}
$\displaystyle\sum_{n\le x} {r_{2,N}(n)}^2$.
\end{center}

A positive squarefree integer $N$ is called solvable (or more
classically ``numerus idoneus'') if $x^2+Ny^2$ has one form per
genus. Note that this means the class number of the form class
group of discriminant $-4N$ equals the number of genera, $2^t$,
where $t$ is the number of distinct prime factors of $N$.
Concerning $r_{2,N}(n)$, Borwein and Choi \cite{BC} proved the
following:

\begin{thm}
Let $N$ be a solvable squarefree integer. Let $x >1$ and $\epsilon
>0$. We have

\begin{center}
$\displaystyle\sum_{n \le x} {r_{2,N}(n)}^2 = \frac{3}{N} \Big(
\displaystyle\prod_{p|2N} \frac{2p}{p+1} \Big ) (x\log x +
\alpha(N)x) + O(N^{\frac{1}{4} + \epsilon} x^{\frac{3}{4} +
\epsilon})$
\end{center}

\noindent where the product is over all primes dividing $2N$ and

\begin{center}
$\alpha(N) = -1 + 2\gamma + \displaystyle\sum_{p|2N} \frac{\log
p}{p+1} + \frac{2L^{\prime}(1,\chi_{-4N})}{L(1, \chi_{-4N})} -
\frac{12}{{\pi}^2} \zeta^{\prime}(2)$
\end{center}

\noindent where $\gamma$ is the Euler-Mascheroni constant and
$L(1, \chi_{-4N})$ is the L-function corresponding to the
quadratic character mod $-4N$.

\end{thm}

Based on this result, Borwein and Choi posed the following:

\begin{conj} For any squarefree $N$,

\begin{center}
$\displaystyle\sum_{n \le x} {r_{2,N}(n)}^2 \sim \frac{3}{N} \Big(
\displaystyle\prod_{p|2N} \frac{2p}{p+1} \Big ) x\log x$
\end{center}

\end{conj}

The main result in \cite{mo} was the following.

\begin{thm}
Let $Q(x,y)=x^2+Ny^2$ for a squarefree integer $N$ with $-N
\not\equiv 1 \bmod 4$. Let $r_{2,N}(n)$ denote the number of
solutions to $n=Q(x,y)$ (counting signs and order). Then
\begin{center}
$\displaystyle\sum_{n \le x} {r_{2,N}(n)}^2 \sim \frac{3}{N}
\Big(\displaystyle \prod_{p|2N} \frac{2p}{p+1} \Big ) x\log x$.
\end{center}
\end{thm}

In this note, we settle the conjecture in the remaining case,
namely

\begin{thm} For $-N \equiv 1 \bmod 4$, we have
\begin{center}
$\displaystyle\sum_{n \le x} {r_{2,N}(n)}^2 \sim \frac{3}{N}
\Big(\displaystyle \prod_{p|2N} \frac{2p}{p+1} \Big ) x\log x$.
\end{center}
\end{thm}

\section{Preliminaries}

Let $Q(x,y)=ax^2 + bxy + cy^2$ denote a positive definite integral
quadratic form with discriminant $D=b^2 - 4ac$ and $gcd(a,b,c)=1$.
Given $Q$, let $\kappa$ be the largest positive integer with
$D/\kappa^2$ an integer congruent to $0$ or $1$ modulo 4. We call
$\kappa$ the {\it conductor} of $Q$ and set $d=D/\kappa^2$. Let
$r(Q,n)$ be the number of representations of the integer $n$ by
the form $Q$. We now relate $r(Q,n)$ to counting the number of
integral ideals of norm $n$ in a given class in a generalized
ideal class group.

Given $D={{\kappa}^2}d$ we consider ideals in $\mathcal{O}_K$
where $K=\mathbb Q(\sqrt{d})$. Let $I_{\kappa}$ be the group of
fractional ideals of $\mathcal{O}_K$ which are quotients of ideals
coprime to $\kappa$ and $P_{\kappa}$ be the subgroup of fractional
ideals which are quotients of principal ideals $\langle \alpha
\rangle \in I_{\kappa}$ where $\alpha \in \mathbb Z +
\kappa\mathcal{O}$. Then set $CL_{\kappa}(K)= I_{\kappa} \diagup
P_{\kappa}$. The elements of $CL_{\kappa}(K)$ correspond
bijectively to proper equivalence classes of positive definite
quadratic forms of discriminant $D={{\kappa}^2}d$. If the proper
equivalence class of $Q$ corresponds to the ideal class
$\frak{c}$, then by \cite{cp}, page 219, we have

\begin{center}
$r(Q,n)= \displaystyle \sum_{r|{\kappa}} w({(\kappa/r)^2}d)
J(\frak{c}_r, n/r^2)$
\end{center}

\noindent where

\begin{center}
$w(D) =\left \{ \begin{array}{l}
6 \quad \mbox{if $D=-3$} \\
4 \quad \mbox{if $D=-4$} \\
2 \quad \mbox{otherwise.}
\end{array}
\right.$ \\
\end{center}

Also $J(\frak{c}_r, n)$ is the number of integral ideals of norm
$n$ in the class $\frak{c}_r$ where $\frak{c}_r$ is the image of
$\frak{c}$ under the natural homomorphism $CL_{\kappa}(K) \to
CL_{{\kappa}/r}(K)$. For the form $Q(x,y)=x^2+Ny^2$ where $-N
\equiv 1 \bmod 4$, the conductor $\kappa=2$ and so we have

$$
\begin{aligned}
r_{2,N}(n) &= w(-4N)J(\frak{c},n) + w(-N)J(\frak{c}_2, n/4) \\
& =2J(\frak{c},n) + w(-N)J(\frak{c}_2, n/4)
\end{aligned}
$$

\noindent where $\frak{c}_2$ is the image under $CL_{2}(K) \to
CL_{1}(K)$, that is, $\frak{c}_2$ is a class in the ideal class
group of $K=\mathbb Q(\sqrt{-N})$.

We now discuss a classical result of Rankin \cite{Rankin} and
Selberg \cite{Sel1} which estimates the size of Fourier
coefficients of a modular form. Specifically, if
$f(z)=\displaystyle\sum_{n=1}^{\infty} a(n) e^{2{\pi}inz}$ is a
nonzero cusp form of weight $k$ on $\Gamma_{0}(N)$, then
\begin{center}
$\displaystyle\sum_{n \leq x} |a(n)|^{2} = \alpha \langle f, f
\rangle x^{k} + O(x^{k-\frac{2}{5}})$
\end{center}

\noindent where $\alpha > 0$ is an absolute constant and $\langle
f, f \rangle$ is the Petersson scalar product. In particular, if
$f$ is a cusp form of weight 1, then $\displaystyle\sum_{n \leq x}
|a(n)|^{2} = O(x)$. One can adapt their result to say the
following. Given two cusp forms of weight $k$ on a suitable
congruence subgroup of $\Gamma=SL_{2}(\mathbb Z)$, say
$f(z)=\displaystyle\sum_{n=1}^{\infty} a(n) e^{2{\pi}inz}$ and
$g(z)=\displaystyle\sum_{n=1}^{\infty} b(n) e^{2{\pi}inz}$, then
\begin{center}
$\displaystyle\sum_{n \leq x} a(n)\overline{b(n)} n^{1-k} = Ax +
O(x^{\frac{3}{5}})$
\end{center}

\noindent where $A$ is a constant. In particular, if $f$ and $g$
are cusp forms of weight 1, then $\displaystyle\sum_{n \leq x}
a(n)\overline{b(n)} = O(x)$.

We conclude this section with a relationship between genus
characters of generalized ideal class groups and the poles of the
Rankin-Selberg convolution of L-functions. Recall that a group
homomorphism $\chi: I_{2} \to S^{1}$ is an ideal class character
if it is trivial on $P_{2}$, i.e.

\begin{center}
$\chi(\langle a \rangle )=1$
\end{center}

\noindent for $a \equiv 1 \bmod \langle 2 \rangle$. Thus an ideal
class character is a character on the generalized class group
$I_{2} \diagup P_{2}$. Recall also that a genus character (see
Chapter 12, section 5 in \cite{iwaniec}) is an ideal class
character of order at most two.

Let us also recall the notion of the Rankin-Selberg convolution of
two L-functions. For squarefree $N$, consider two ideal class
characters $\chi_1$, $\chi_2$ for $CL_{2}(K)$, the generalized
ideal class group of $K=\mathbb Q(\sqrt{-N})$ and their associated
Hecke L-series

\begin{center}
$L_{2}(s, \chi_1)= \displaystyle \sum_{(\frak{a}, 2)=1}
\frac{\chi_1(\frak{a})}{{N(\frak{a})}^s}$
\end{center}

\begin{center}
$L_{2}(s, \chi_2)= \displaystyle \sum_{(\frak{a}, 2)=1}
\frac{\chi_2(\frak{a})}{{N(\frak{a})}^s}$
\end{center}

\noindent which converge absolutely in some right half-plane. We
form the convolution L-series by multiplying the coefficients,

\begin{center}
$L_{2}(s, \chi_1 \otimes \chi_2) = \displaystyle \sum_{(\frak{a},
2)=1} \frac{\chi_1(\frak{a})\chi_2(\frak{a})}{{N(\frak{a})}^s}$
\end{center}

The following result describes a relationship between genus
characters $\chi$ and the orders of poles of $L_{2}(s, \chi
\otimes \chi)$. The proof is similar to that of Proposition 2.4 in
\cite{mo}.

\begin{pro} Let $\chi$ be an ideal class character for $CL_{2}(K)$,
$-N \equiv 1 \bmod 4$, and $L_{2}(s,\chi)$ the associated Hecke
L-series. Then $\chi$ is a genus character if and only if
$L_{2}(s, \chi \otimes \chi)$ has a double pole at $s=1$.
\end{pro}

\begin{rem}
By Proposition 2.1, if $\chi$ is a non-genus character, then
$L_{2}(s, \chi \otimes \chi)$ has at most a simple pole at $s=1$.
\end{rem}

\section{Proof of Theorem 1.4}

\begin{proof}
As the proof is similar to that of Theorem 1.3 in \cite{mo}, we
sketch the relevant details. If $-N \equiv 1 \bmod 4$, then the
discriminant of $K=\mathbb Q(\sqrt{-N})$ is $-N$. We also assume
that $t$ is the number of distinct prime factors of $N$ and so the
discriminant $-N$ also has $t$ distinct prime factors. For
$K=\mathbb Q(\sqrt{-N})$, consider the zeta function

\begin{center}
$\zeta_{K}(s, 2) = \displaystyle\sum_{(\frak{a}, 2)=1}
\frac{1}{{N(\frak{a})}^s}$
\end{center}

\noindent where the sum is over those ideals $\frak{a}$ of
$\mathcal{O}_{K}$ prime to $2$. We now split up $\zeta_{K}(s, 2)$,
according to the classes $\frak{c}_{i}$ of the generalized ideal
class group $CL_{2}(K)$, into the partial zeta functions (see page
161 of \cite{lang})
\begin{center}
$\zeta_{\frak{c}_i}(s) = \displaystyle\sum_{\frak{a} \in
\frak{c}_{i}} \frac{1}{{N(\frak{a})}^s}$
\end{center}

\noindent so that $\zeta_{K}(s,2)=
\displaystyle\sum_{i=0}^{h_{2}-1} \zeta_{\frak{c}_i}(s)$ where
$h_{2}$ is the order of $CL_{2}(K)$.

Let $\frak{c}$ be the ideal class in $CL_{2}(K)$ which corresponds
to the proper equivalence class of $Q(x,y)=x^2+Ny^2$. Now let
$\chi$ be an ideal class character of $CL_{2}(K)$ and consider the
Hecke L-series for $\chi$, namely

\begin{center}
$L_{2}(s, \chi)= \displaystyle \sum_{(\frak{a}, 2)=1}
\frac{\chi(\frak{a})}{{N(\frak{a})}^s}$.
\end{center}

\noindent We may now rewrite the Hecke L-series as

\begin{center}
$L_{2}(s,\chi) = \displaystyle\sum_{i=0}^{h_{2}-1}
\chi(\frak{c}_i) \zeta_{\frak{c}_{i}}(s)$.
\end{center}

\noindent And so summing over all ideal class characters of
$CL_{2}(K)$, we have

\begin{center}
$\displaystyle\sum_{\chi} \overline{\chi}(\frak{c}) L_{2}(s, \chi)
= \displaystyle\sum_{i=0}^{h_{2}-1} \zeta_{\frak{c}_i}(s) \Big (
\displaystyle\sum_{\chi} \overline{\chi}(\frak{c})
\chi(\frak{c}_{i}) \Big )$.
\end{center}

\noindent The inner sum is nonzero precisely when
$\frak{c}=\frak{c}_{i}$. Thus we have

\begin{center}
$\zeta_{\frak{c}}(s) = \displaystyle \frac{1}{h_{2}} \sum_{\chi}
\overline{\chi}(\frak{c}) L_{2}(s, \chi)$
\end{center}

\noindent and so

\begin{center}
$\zeta_{\frak{c}}(s) = \displaystyle \frac{1}{h_2} (L_{2}(s,
\chi_{0}) + \overline{\chi_1}(\frak{c})L_{2}(s, \chi_{1}) + \dots
+ \overline{\chi_{h_{2}-1}}(\frak{c})L_{2}(s, \chi_{{h_{2}}-1}))$.
\end{center}

As $\chi_{0}$ is the trivial character, $L_{2}(s, \chi_{0}) =
\zeta_{K}(s,2)$. Comparing $n^{th}$ coefficients, we have

\begin{center}
$J(\frak{c},n)= \displaystyle \frac{1}{h_{2}} (a_{n} + b_{1}(n) +
\dots + b_{h_{2}-1}(n))$.
\end{center}

\noindent where $a_{n}$ is the number of integral ideals of
$\mathcal{O}_{K}$ prime to 2 and of norm $n$ and the $b_{i}$'s are
coefficients of weight 1 cusp forms (see \cite{bump}). Recall we
also have

\begin{center}
$r_{2,N}(n)= 2J(\frak{c},n) + w(-N)J(\frak{c}_2, n/4)$
\end{center}

\noindent and so

\begin{center}
$r_{2,N}(n)=\displaystyle \frac{2}{h_2}\Big( a_{n} + b_{1}(n) +
\dots + b_{h_{2}-1}(n)\Big ) + w(-N)J(\frak{c}_2, n/4)$.
\end{center}

Thus

\begin{center}
$\displaystyle\sum_{n \le x} {r_{2,N}(n)}^2=\frac{4}{h_{2}^2}
\Big( \sum_{n\le x} a_{n}^2 + \sum_{\substack{i \\ n\le x }}
{b_{i}(n)}^2 + 2 \sum_{\substack{i \\ n\le x }} a_{n}b_{i}(n) +
\sum_{\substack{i \neq j \\ n\le x}} b_{i}(n)b_{j}(n) \Big) +
\frac{4}{h_{2}} \sum_{n \le x}\Big ( a_{n} + b_{1}(n) + \dots +
b_{h_{2}-1}(n) \Big )w(-N)J(\frak{c}_2, n/4) + \sum_{n \le
x}w(-N)^2 J(\frak{c}_2, n/4)^2$.
\end{center}

Assume $-N \equiv 1 \bmod 8$. Applying the main theorem in
\cite{Nowak} to the Dirichlet series $\displaystyle
\sum_{n=1}^{\infty} \frac{a_{n}^2}{n^s}$, we obtain

\begin{center}
$\displaystyle \sum_{n\leq x} a_{n}^2 \sim Ax\log x$
\end{center}

\noindent where $A=\displaystyle \frac{1}{2\pi^2}{L(1,
\chi_{-N})}^2 \prod_{p|N} \frac{p}{p+1}$. As $-N$ has $t$ distinct
prime factors, we have $2^t$ genus characters for $CL(K)$ where
$K=\mathbb Q(\sqrt{-N})$. By \cite{lang} (see Theorem 1, page
127), we have $2^t$ genus characters for $CL_2(K)$. We now must
estimate $\displaystyle \sum_{\substack{i \\ n\le x }}
{b_{i}(n)}^2$. Let us now assume that the first $2^{t} - 1$ terms
arise from L-functions associated to genus characters. By
Proposition 2.1 and an application of Perron's formula, we obtain

\begin{center}
$\displaystyle \sum_{n\le x} {b_{i}(n)}^2 \sim Ax\log x$.
\end{center}

\noindent As this estimate holds for each $i$ such that $1 \le i
\le 2^{t}-1$, the term $Ax \log x$ appears $2^{t}$ times in the
estimate of $\displaystyle\sum_{n \le x} {r_{2,N}(n)}^2$. By
Remark 2.2 and the Rankin-Selberg estimate, the remaining terms
are all $O(x)$. Thus

\begin{center}
$\displaystyle\sum_{n \le x} {r_{2,N}(n)}^2 \sim \frac{4}{h_{2}^2}
\Big ( 2^{t} \frac{1}{2{\pi}^{2}} {L(1, \chi_{-N})}^2
\displaystyle\prod_{p|N} \frac{p}{p+1} \Big ) x \log x$.
\end{center}

By \cite{cohn}, we have $L(1, \chi_{-N})=\frac{h\pi}{\sqrt{N}}$
where $h$ is the class number of $K$ and $h_2=h$. Thus

\begin{center}
$\displaystyle\sum_{n \le x} {r_{2,N}(n)}^2 \sim \frac{3}{N}
\Big(\displaystyle \prod_{p|2N} \frac{2p}{p+1} \Big ) x\log x$.
\end{center}

For $-N \equiv 5 \bmod 8$, we have $h_2=3h$ and again by
\cite{Nowak},

\begin{center}
$\displaystyle \sum_{n\leq x} a_{n}^2 \sim \Big(
\frac{9}{2\pi^2}{L(1, \chi_{-N})}^2 \prod_{p|N} \frac{p}{p+1}\Big)
x\log x$.
\end{center}

Thus

\begin{center}
$\displaystyle\sum_{n \le x} {r_{2,N}(n)}^2 \sim \frac{3}{N}
\Big(\displaystyle \prod_{p|2N} \frac{2p}{p+1} \Big ) x\log x$.
\end{center}

\end{proof}

\begin{rem}
We would like to mention another approach which confirms Theorems
1.3 and 1.4. Let $Q \in \mathbb Z^{2 \times 2}$ be a non-singular
symmetric matrix with even diagonal entries and $q({\bf
x})=\frac{1}{2}Q[{\bf x}]=\frac{1}{2}{{ \bf x}^T}Q{\bf x}$, ${\bf
x} \in \mathbb Z^{2}$, the associated quadratic form in two
variables. Let $r(Q,n)$ denote the number of representations of
$n$ by the quadratic form $Q$. Now consider the theta function

\begin{center}
$\theta_{Q}(z)=\displaystyle \sum_{{\bf x} \in \mathbb Z^{2}}
e^{{\pi}izQ[{\bf x}]}$.
\end{center}

\noindent The Dirichlet series associated with the automorphic
form $\theta_{Q}$ is

\begin{center}
$(4\pi)^{-1/2}\zeta_{Q}(\frac{1}{2}+ s)$
\end{center}

\noindent where

\begin{center}
$\zeta_{Q}(s)=\displaystyle \sum_{n=1}^{\infty} \frac{r(Q,n)}{n^s}
= \sum_{{\bf x} \in \mathbb Z^{2}\setminus\{\bf{0}\}} q({\bf
x})^{-s}$
\end{center}

\noindent for $\Re(s)>1$. A careful and involved application of
the Rankin-Selberg method to the above Dirichlet series (see
Theorems 2.1 and 5.1 in \cite{mul1} and Theorem 5.2 in
\cite{mul2}) combined with a Tauberian argument yields the
following (see Theorem 6.1 in \cite{mul1})

\begin{center}
$\displaystyle \sum_{n\leq x} r(Q,n)^2 \sim A_{Q}x\log x$
\end{center}

\noindent where

\begin{center}
$\displaystyle A_{Q}=12\frac{A(q)}{q} \prod_{p \mid q} \Big(1+
\frac{1}{p}\Big )^{-1}$.
\end{center}

\noindent Here $q=$ det $Q$ and $A(q)$ denotes the multiplicative
function defined by

\begin{center}
$A(p^e)=2 + (1-\frac{1}{p})(e-1)$
\end{center}

\noindent where $p$ is an odd prime, $e\geq 1$, and

$$
A(2^e) =\left \{ \begin{array}{l}
1 \quad \mbox{if $e \leq 1$},\\
2 \quad  \mbox{if $e=2$},\\
e-1 \quad \mbox{if $e \geq 3$.}
\end{array}
\right. \\
$$

Let us now turn to our situation. Consider $q({\bf
x})=x^2+Ny^2=\frac{1}{2}{{\bf x}^T}Q{\bf x}$ where
$Q=\left( \begin{matrix} 2 & 0 \\
0 & 2N \\
\end{matrix} \right)$,
$N$ squarefree. Thus $q=4N$. Suppose $N$ has $t$ distinct prime
factors. Then $A(4N)=2^{t+1}$ and so

\begin{center}
$A_{Q}=\displaystyle \frac{3}{N}2^{t+1}\prod_{p\mid
2N}\Big(1+\frac{1}{p}\Big)^{-1}=\frac{3}{N}\prod_{p\mid 2N}
\frac{2p}{p+1}$.
\end{center}

\end{rem}

\section*{Acknowledgments}
The author would like to thank Ram Murty for many productive
discussions and for his comments on a preliminary version of this
paper.


\begin{thebibliography}{10}
\bibitem{BC}
J. Borwein, K.K. Choi, \emph{On Dirichlet Series for sums of
squares}, Rankin memorial issues, Ramanujan J. \textbf{7} (2003),
no.1-3, 95--127.

\bibitem{bump}
D. Bump, {\em Automorphic forms and representations}, Cambridge
Studies in Advanced Mathematics, \textbf{55}, Cambridge University
Press, 1997.

\bibitem{cp}
R. Chapman, A. van der Poorten, \emph{Binary Quadratic Forms and
the Eta Function}, Number theory for the millennium, I (Urbana,
IL, 2000), 215--227, A K Peters, Natick, MA, 2002.

\bibitem{cohn}
H. Cohn, {\em Advanced Number Theory}, Dover Publications, Inc.,
New York, 1980.

\bibitem{iwaniec}
H. Iwaniec, {\em Topics in Classical Automorphic Forms}, Graduate
Studies in Mathematics, Vol. 17, Amer. Math. Soc., Providence, RI,
1997.

\bibitem{Nowak}
M. K$\ddot{u}$hleitner, W.G. Nowak, \emph{The average number of
solutions to the Diophantine equation $U^2+V^2=W^3$ and related
arithmetic functions}, Acta Math. Hungar. \textbf{104} 
(2004), 225--240.

\bibitem{lang}
S. Lang, {\em Algebraic Number Theory}, Second Edition,
Springer-Verlag, New York, 1994.

\bibitem{mul1}
W. M$\ddot{u}$ller, \emph{The mean square of Dirichlet series
associated with automorphic forms}, Monatsh. Math. \textbf{113}
(1992), 121--159.

\bibitem{mul2}
W. M$\ddot{u}$ller, \emph{The Rankin-Selberg Method for
non-holmorphic automorphic forms}, J. Number Theory \textbf{51}
(1995), 48--86.

\bibitem{mo}
R. Murty, R. Osburn, \emph{Representations of integers by certain
positive definite binary quadratic forms}, submitted.

\bibitem{Rankin}
R.A. Rankin, \emph{Contributions to the theory of Ramanujan's
function $\tau(n)$ and similar functions. II. The order of the
Fourier coefficients of integral modular forms}, Proc. Cambridge
Philos. Soc. \textbf{35} (1939), 357--373.

\bibitem{Sel1}
A. Selberg, \emph{Bemerkungen $\ddot{u}$ber eine Dirichletsche
Reihe, die mit der Theorie der Modulformen nahe verbunden ist},
Archiv. Math. Natur. B \textbf{43} (1940), 47--50.

\end{thebibliography}
\end{document}